\documentclass[leqno]{amsart}

\usepackage{algpseudocode}
\usepackage{amssymb}
\usepackage{amsmath}

\def\bne{\begin{equation}}
\def\ene{\end{equation}}
\def\F{\mathbb{F}}
\def\Re{\mathbb{R}}
\DeclareMathOperator{\ulp}{ulp}
\DeclareMathOperator{\bu}{\mathbf{u}}

\newtheorem{algorithm}{Algorithm}

\title{Fast Compensated Algorithms for the Reciprocal Square Root, the Reciprocal Hypotenuse, and Givens Rotations}

\author{Carlos F. Borges\\Department of Applied Mathematics\\Naval Postgraduate School\\Monterey CA 93943\\borges@nps.edu}

\begin{document}

\begin{abstract}
The reciprocal square root is an important computation for which many very sophisticated algorithms exist (see for example \cite{Moroz,863046,863031} and the references therein). In this paper we develop a simple differential compensation (similar to those developed in \cite{borges}) that can be used to improve the accuracy of a naive calculation. The approach relies on the use of the fused multiply-add (FMA) which is widely available in hardware on a variety of modern computer architectures. We first show how compensate by computing an exact Newton step and investigate the properties of this approach.  We then show how to leverage the exact Newton step to get a modified compensation which requires one additional FMA and one additional multiplication. This modified method appears to give correctly rounded results experimentally and we show that it can be combined with a square root free method for estimating the reciprocal square root to get a method that is both very fast (in computing environments with a slow square root) and, experimentally, highly accurate. Finally, we show how these approaches can be extended to the reciprocal hypotenuse calculation and the construction of Givens rotations.
\end{abstract}

\maketitle

There are many current algorithms (see for example \cite{Moroz,863046,863031} and the references therein) in the literature which can be used to compute the reciprocal square root. In this paper we will show how one can leverage the fused multiply-add to compute a differential compensation, that turns out to be an exact Newton step, that will yield a very accurate answer. We then show how to leverage the exact Newton step to get a modified compensation which requires one additional FMA operation. When this approach is combined with a square root free method for estimating the reciprocal square root we get a method that is both very fast and, experimentally, highly accurate. Such methods are important in computing environments that do not have a fast square root (e.g. microcontrollers and FPGAs). We then show how this same device can be used to fast and accurate algorithms for computing the reciprocal hypotenuse and constructing Givens rotations. We give a careful error analysis only for the reciprocal square root calculation, but we illustrate the accuracy of all of these algorithms {\em experimentally} by comparing them to extended precision calculations carried out with the MPFR package in Julia. All of these algorithms can be implemented on architectures that do not have a hardware FMA by using a software implementation of the FMA although this would likely signficantly reduce the speed.

This paper will be restricted to the case where all floating-point calculations are done in IEEE 754 compliant radix~2 arithmetic using round-to-nearest rounding, although many of the results can be extended to other formats under proper conditions.

\section{Mathematical Preliminaries}

We begin with a few definitions. We shall denote by $\F \subset \Re$ the set of all radix 2 floating-point numbers with precision $p$. Because all of the algorithms we are investigating can avoid overflow and underflow by judicious scaling with exact powers of 2 we assume that the exponent range is infinite. We define $fl(x):\Re \rightarrow \F$ to be a function such that $fl(x)$ is the element of $\F$ that is closest to $x$. Since we have restricted ourselves to radix 2 IEEE754 compliant arithmetic we will assume that round-to-even is used in the event of a tie. Throughout this paper we assume that $ulp(x)$ is a unit in the last place as defined in \cite{Goldberg}. To wit, we define $ulp(x) = 2^{e-p+1}$ for any $x \in [2^e,2^{e+1})$ where $p$ is the precision of the floating-point format. For example, in double precision $\ulp(1) = 2^{-52}$. We define the {\em unit roundoff}, which we will denote by $\bu = \frac{1}{2}\ulp(1)$ so that, for example, in double precision $\bu = 2^{-53}$.

\section{Computing the reciprocal of the square root}

Given a floating-point number $x \in \F$ we wish to compute $$y = \frac{1}{\sqrt{x}}$$ in floating point. The most accurate naive approach is simply\footnote{One can also use y=1/sqrt(x) but it can lead to errors greater than 1 ulp. See \cite{Markstein2}.}

\vspace{.2in}
\begin{algorithm} Naive rsqrt

\hrule
\begin{algorithmic}
	\State {r = 1/x}
	\State {y = sqrt(r)}
\end{algorithmic}
\hrule
\label{Naiversqrt}
\end{algorithm}
\vspace{.2in}

We note that the computed quantity $y$ is subject to various errors due to the effects of finite precision and hence
\bne
\sqrt{\frac{1}{x}} = y (1+\nu)
\label{recip1}
\ene
for some $\nu \in \Re$, where it is reasonable to assume that $|\nu| < 1$. Under our assumptions as to the floating-point environment, a standard error analyis reveals that 
\bne
|\nu| < \frac{3}{2}\bu + O(\bu^2).
\label{crap}
\ene

Squaring both sides of \ref{recip1} and a bit of algebra gives us
\bne
x y^2 = \frac{1}{(1+\nu)^2}
\label{recip2}
\ene
and then one step of long division on the right hand side and a bit more algebra yields
\bne
\nu = \frac{1-x y^2}{2}  + \nu^2\frac{3-2\nu}{2(1+\nu)^2}.
\label{recip3}
\ene

Observe that at this point that for any $y$ satisfying \ref{recip1} with $|\nu| < 1$ we have
\bne
\sqrt{\frac{1}{x}} = y + y\left(\frac{1-xy^2}{2} + \nu^2\frac{3-2\nu}{2(1+\nu)^2}\right).
\label{bigone}
\ene
If we ignore the $O(\nu^2)$ term and let
\bne
\bar{\nu} = \frac{1-xy^2}{2}.
\label{nu}
\ene
We can add a compensation to the naively computed value using
\bne
y_C = y + y\bar{\nu}
\label{Newtonstep}
\ene
which the reader may recognize as the Newton iteration for
$$
f(y) = x - \frac{1}{y^2}
$$
which is a common approach to estimating the reciprocal square root (see \cite{Moroz} and references therein).

We can perform a traditional rounding error analysis for the single Newton step in \ref{Newtonstep} using the standard model for radix 2 round-to-nearest floating point arithmetic. Going forward note that all $|\sigma_i| \leq \bu$. Assume that the computed value of $\bar{\nu}$ is  $\bar{\nu}(1+\epsilon)$ for some $|\epsilon| < 1$. Then, if we apply the Newton step using an FMA the final computed value, which we shall call $y_F$, satisfies
\begin{eqnarray*}
y_F & = &  (y + y  \frac{1-xy^2}{2}(1+\epsilon))(1+\delta_1) \\
& = &  y + y \frac{1-xy^2}{2}  + y\delta_1 + y \frac{1-xy^2}{2}(\epsilon (1+ \delta_1)). \\
\end{eqnarray*}
Using \ref{bigone} to replace the first two terms on the right we get
\bne
y_F =  \sqrt{\frac{1}{x}} -y\nu^2\frac{3-2\nu}{2(1+\nu)^2}  + y\delta_1 + y \frac{1-xy^2}{2}(\epsilon (1+ \delta_1))
\ene
and then using \ref{recip3} gives
\begin{eqnarray*}
y_F & = &  \sqrt{\frac{1}{x}} -y\nu^2\frac{3-2\nu}{2(1+\nu)^2}  + y\delta_1 + y \left(\nu - \nu^2\frac{3-2\nu}{2(1+\nu)^2}\right)(\epsilon (1+ \delta_1)) \\
& = &  \sqrt{\frac{1}{x}} + y \left( \delta_1 - \nu^2\frac{3-2\nu}{2(1+\nu)^2} + \left(\nu - \nu^2\frac{3-2\nu}{2(1+\nu)^2}\right)(\epsilon (1+ \delta_1))\right)\\
& = &  \sqrt{\frac{1}{x}} + y \left( \delta_1 + O(\nu^2) + O(\nu \epsilon) \right)
\end{eqnarray*}
And finally
$$
\frac{\left| y_F - \sqrt{\frac{1}{x}}\right|}{\sqrt{\frac{1}{x}}}  = \frac{\left|y\left( \delta_1 + O(\nu^2) + O(\nu \epsilon) \right)\right|}{y(1+\nu)} = \frac{\left|\left( \delta_1 + O(\nu^2) + O(\nu \epsilon) \right)\right|}{(1+\nu)}.
$$
Clearly, if both $\nu$ and $\epsilon$ are $O(\bu)$ then 
\bne
y_F = \sqrt{\frac{1}{x}} + y \left( \delta_1 + O(\bu^2)\right).
\label{goodbound}
\ene
This is a useful condition as it implies that any $y_F$ that is not correctly rounded is directly adjacent to the correctly rounded value\footnote{A somewhat stronger statement than {\em faithful}.} and, further, that the true value lies very near the center of the circumscribed interval. Whenever $y_F$ satisfies \ref{goodbound} we shall say that it is {\em weakly rounded}. However, the problem is that when $\nu$ gets very small, the relative error, $\epsilon$ in computing $\bar{\nu}$ begins to climb as this computation is subject to extreme cancellation. To avoid this disaster we rewrite the subject quantity in the following suggestive form
$$
\bar{\nu} = \frac{1}{2}(1-xr - x(y^2-r))
$$
and propose the following compensated algorithm which requires adding four FMA calls and a single multiply to the naive algorithm:

\vspace{.2in}
\begin{algorithm} Compensated rsqrt
\hrule
\begin{algorithmic}
   \State {r = 1/x}
	\State {y = sqrt(r)}
	\State {mxhalf = -0.5*x}
   \State {$\sigma$ = fma(mxhalf, r, 0.5)}
   \State {$\tau$ = fma(y, y, -r)}
   \State {$\bar{\nu}$ = fma(mxhalf, $\tau$, $\sigma$)}
   \State {y = fma(y, $\bar{\nu}$, y)}
\end{algorithmic}
\hrule
\label{Correctedrsqrt}
\end{algorithm}
\vspace{.2in}

\subsection{Error Analysis}

To see why algorithm \ref{Correctedrsqrt} works we  simply need to observe that under appropriate conditions (see Chapter 4 Theorems 4.9 and 4.10 in \cite{FPHB}), both $\sigma, \tau \in \F$,  and are computed exactly using the FMA (see \cite{Markstein} or Chapter 4 Corollary 4.11 and 4.12 in \cite{FPHB} for the precise conditions). This means that the computed value of $\bar{\nu}$ is correctly rounded (using an FMA) and hence $\epsilon = \bu$ which implies that the result of algorithm \ref{Correctedrsqrt} satisfies the bound in equation \ref{goodbound}.

\subsection{Numerical Testing}

Since our error bound is not sufficiently tight to show correct rounding, we now test the compensated algorithm against the naive approach. All testing is done in IEEE754 double precision arithmetic with code written in Julia 1.5.3 running on an Intel(R) Core(TM) i7-7700K CPU @ 4.20GHz. As a baseline for testing purposes we will compute the reciprocal square root, $\bar{y}$, by using the BigFloat format in Julia which uses the GNU MPFR package to do an extended precision calculation. It is critical to note that neither square root nor division are finite operations. That is, both can result in infinite length results and hence this extended precision computation is necessarily prone to double rounding. That means that the result of the MPFR computation cannot be guaranteed to represent a correctly rounded value of the true result (although it will do so in nearly every case).

 We will do this using $10^{9}$ uniformly distributed double precision (Float64) random inputs. We will run both algorithms on each random input as well as computing $\bar{y}$. In the tables below we summarize the error rates of each algorithm which we define to be the percentage of times each algorithm differed from $\bar{y}$ by exactly zero ulp and exactly one ulp. We note that neither algorithm ever differed from $\bar{y}$ by more than one ulp.

\begin{table}[h]
  \begin{center}
    \caption{Error rate for computing the reciprocal square root with $x \sim \mathcal{U}(1/2,1)$}
    \label{tab:tablersqrt1}
    \begin{tabular}{l|l|l}
  & Naive & Compensated \\
      \hline
Zero ulp & 89.227  &  100\\
One ulp & 10.773 & 0
    \end{tabular}
  \end{center}
\end{table}

\begin{table}[h]
  \begin{center}
    \caption{Error rate for computing the reciprocal square root with $x \sim \mathcal{U}(1,2)$}
    \label{tab:tablersqrt2}
    \begin{tabular}{l|l|l}
  & Naive & Compensated \\
      \hline
Zero ulp & 84.762  &  100\\
One ulp  & 15.238 & 0
    \end{tabular}
  \end{center}
\end{table}

Although the compensated algorithm is clearly very accurate, it does not always yield a correctly rounded result. To wit, let $x = 1-2\bu$, then the computed value of $y=1.0$ and the compensation will be computed exactly as $y\bar{\nu} = \bu$, and $1+\bu \rightarrow 1$ because of round-to-even.\footnote{We are not aware of any other examples beyond this one and any multiples of the form $(1-2\bu)4^k$.} Since we know that the reciprocal square root of any radix 2 floating point number cannot be the exact midpoint of two consecutive floating point numbers (see Chapter 4 Theorem 4.20 in \cite{FPHB}) it is clear that the problem in this case is the difference between $\bar{\nu}$ and $\nu$.

Note that  $\bar{\nu}$ and $\nu$ always have the same sign, and further that
\bne
\bar{\nu} < \nu
\label{overcorrect}
\ene
provided that $\nu \neq 0$,which is seen by rearranging \ref{recip3}. This means that algorithm \ref{Correctedrsqrt} will always slightly undercompensate if $\bar{\nu}$ is positive, and slightly overercompensate if it is negative. Now, in the specific example just discussed, the compensation was positive and hence too small and we really should have compensated by a tad more which would have given us the correctly rounded result of $y=1+2\bu$. We note that although the result generated by the algorithm is not correctly rounded, it is {\em weakly rounded} as guaranteed by the earlier analysis.

\subsection{A Modified Compensation}

The fact that $\bar{\nu}$ can be computed exactly means it may reasonably be used to estimate $\nu$ directly. To that end, if $\nu = O(\bu)$ then equation \ref{recip3} can be rewritten as
$$
\nu = \bar{\nu}  + \frac{3}{2}\nu^2 + O(\bu^3).
$$
A little algebra introduces some more $O(\bu^3)$ terms and yields
\begin{eqnarray*}
\nu  & = & \frac{\bar{\nu}}{1 - \frac{3}{2}\nu} + O(\bu^3)\\
& = & \bar{\nu} \left(1 + \frac{3}{2}\nu \right) + O(\bu^3)
\end{eqnarray*}
and finally, replacing $\nu$ on the right hand side with $\bar{\nu}$ introduces yet more $O(\bu^3)$ terms and gives
\bne
\nu  = \bar{\nu} \left(1 + \frac{3}{2}\bar{\nu} \right) + O(\bu^3).
\label{ohyeah}
\ene

This gives a modified compensation form of the algorithm which the careful reader will notice is simply Halley's method.

\vspace{.2in}
\begin{algorithm} Modified Compensated rsqrt
\hrule
\begin{algorithmic}
   \State {r = 1/x}
	\State {y = sqrt(r)}
	\State {mxhalf = -0.5*x}
   \State {$\sigma$ = fma(mxhalf, r, 0.5)}
   \State {$\tau$ = fma(y, y, -r)}
   \State {$\bar{\nu}$ = fma(mxhalf, $\tau$, $\sigma$)}
    \State {$\nu$ = fma(1.5*$\bar{\nu}$,$\bar{\nu}$,$\bar{\nu}$)}
   \State {y = fma(y, $\nu$, y)}
\end{algorithmic}
\hrule
\label{Correctedrsqrt2}
\end{algorithm}
\vspace{.2in}

This algorithm requires one more multiply and one more FMA than the original. In testing this algorithm we have never found a single example where it fails to generate the correctly rounded result, it even does so for $x = 1-2\bu$. The author conjectures that this algorithm always gives a correctly rounded result but does not have a proof.

\subsection{Developing a fast square-root free variant}

The compensation can be applied to existing square-root-free algorithms for computing the reciprocal square root. Such methods are important in computing environments that do not have a fast square root (e.g. microcontrollers and FPGAs). There are many such algorithms but we have chosen to use one of a type that is best known for appearing in the code for the video game {\em Quake III Arena}. The specific algorithm from \cite{Moroz} that we will investigate is called {\tt RcpSqrt331d(x)}. A Julia port of that code appears below.
\vspace{.1in}

\begin{verbatim}
function RcpSqrt331d(x::Float64)
    i = reinterpret(Int64,x)
    k = i & 0x0010000000000000
    if k != 0
        i = 0x5fdb3d14170034b6 - (i >> 1)
        y = reinterpret(Float64,i)
        y = 2.33124735553421569*y*fma(-x, y*y,1.07497362654295614 )
    else
        i = 0x5fe33d18a2b9ef5f - (i >> 1)
        y = reinterpret(Float64,i)
        y = 0.82421942523718461*y*fma(-x, y*y, 2.1499494964450325)
    end
    mxhalf = -0.5*x
    y = y*fma(mxhalf, y*y, 1.5000000034937999)
  # The next two lines are a single step of Newton
    r = fma(mxhalf, y*y, 0.5)
    y = fma(y, r, y)
end
\end{verbatim}
\vspace{.1in}

Note that the final two lines of the code are simply a careful application of one step of Newton which requires two FMAs and one multiply. We replace those lines with the modified form of compensation. The added lines are:

\vspace{.2in}
\begin{algorithm} {\tt RcpSqrt331dModified(x) - partial}
\hrule
\begin{algorithmic}
	\State{\# The following lines replace the Newton step in RcpSqrt331d(x)}
   \State {r = 1/x}
   \State {$\sigma$ = fma(r,mxhalf,.5)}
   \State {$\tau$ = fma(y, y, -r)}
   \State {$\bar{\nu}$ = fma(mxhalf, $\tau$, $\sigma$)}
    \State {$\nu$ = fma(1.5*$\bar{\nu}$,$\bar{\nu}$,$\bar{\nu}$)}
   \State {y = fma(y, $\nu$, y)}
\end{algorithmic}
\hrule
\label{FastCorrectedrsqrt}
\end{algorithm}
\vspace{.2in}

This adds three FMAs and one divide to the code and benchmark timings in Julia indicate that it adds about 5\% to the execution time.

In the table below we show the results from the accuracy test comparing the two algorithms. Note that the compensated form returns the correctly rounded answer every single time in the experiment. Moreover, it gives a correctly rounded answer for $x = 1-2\bu$.

\begin{table}[h]
  \begin{center}
    \caption{Error rate for computing the reciprocal square root with $x \sim \mathcal{U}(1/2,1)$}
    \label{tab:tablefrsqrt1}
    \begin{tabular}{l|l|l}
  & RcpSqrt331d & RcpSqrt331dModified \\
      \hline
Zero ulp & 87.324  &  100\\
One ulp & 12.676 & 0
    \end{tabular}
  \end{center}
\end{table}

\begin{table}[h]
  \begin{center}
    \caption{Error rate for computing the reciprocal square root with $x \sim \mathcal{U}(1,2)$}
    \label{tab:tablefrsqrt2}
    \begin{tabular}{l|l|l}
  & RcpSqrt331d & RcpSqrt331dModified \\
      \hline
Zero ulp & 82.119  &  100\\
One ulp  & 17.881 & 0
    \end{tabular}
  \end{center}
\end{table}

\section{Computing the reciprocal of the hypotenuse}

This same device can be used to correct a naive algorithm that computes the reciprocal of the hypotenuse of a right triangle with sides $x$ and $y$ which is given by $$\rho = \frac{1}{\sqrt{x^2+y^2}}.$$ Such a function is generally called {\tt rhypot(x,y)} and it is important to note that the order of the arguments and their signs should not change the mathematical value of the output. However, order is important in our computations ansd so we will enforce the condition that $x \geq y > 0$ within our codes. It is also important to note that it may be necessary to rescale the problem and then scale back the answer to avoid intermediate overflow/underflow issues. This process is carefully discussed in \cite{borges} and we refer the reader to that paper for the details. We do note that rescaling is rarely needed and so it is generally better to write algorithms that trap floating point exceptions and only rescale in the event of an overflow/underflow error in an intermediate step.

The naive approach is to simply compute

\vspace{.2in}
\begin{algorithm} Naive rhypot

\hrule
\begin{algorithmic}
   \State {$r = 1/(x^2+y^2)$}
	\State {$\rho = sqrt(r)$}
\end{algorithmic}
\hrule
\label{Naiverhypot}
\end{algorithm}
\vspace{.2in}

To construct a compensated algorithm we begin by noting that $\tau$, in this case is the same as before. However, the value of $\sigma$ must now satisfy
\bne
r(x^2+y^2) + \sigma = 1.
\label{hypeq}
\ene
It is important to note that we no longer have any assurance that $\sigma$ is a floating point number since $x^2+y^2$ may not be one. And although we cannot generally compute it exactly, it is possible to compute it to high relative accuracy if we do so carefully as in \cite{borges}. Using those techniques leads to a compensated algorithm.

\vspace{.2in}
\begin{algorithm} Compensated rhypot

\hrule
\begin{algorithmic}
	\State {$x = abs(x)$}
	\State {$y = abs(y)$}
	\If{$x<y$}
	\State {$swap(x,y)$}
   \EndIf
	\State {$x_{sq} = x^2$}
	\State {$y_{sq} = y^2$}
	\State {$\sigma = x_{sq}+y_{sq}$}
	\State {$\sigma_e = y_{sq}-(sigma-x_{sq})+fma(x,x,-x_{sq})+fma(y,y,-y_{sq})$}
   \State {$r = 1/\sigma$}
   \State {$\sigma = fma(-r,\sigma_e,fma(-r,\sigma,1))$}
	\State {$\rho = sqrt(r)$}
	\State {$\tau = fma(-\rho,\rho,r)$}
	\State {$\nu = fma(\sigma,\tau,\sigma)/2$}
	\State {$\rho =fma(\rho,\nu,\rho)$}
\end{algorithmic}
\hrule
\label{Correctedrhypot}
\end{algorithm}
\vspace{.2in}

Once again, as a baseline for testing purposes we will first perform the naive calculation using the BigFloat format in Julia to compute $\bar{\rho}$. And once again we note that this computation cannot be guaranteed to represent a correctly rounded value of the true result (although it will do so in nearly every case). We will use $10^{9}$ normally distributed random inputs, that is both $x,y \sim \mathcal{N}(0,1)$. In the table below we summarize the error rates of each algorithm.

\begin{table}[h]
  \begin{center}
    \caption{Error rate for computing the reciprocal hypotenuse with $x,y \sim \mathcal{N}(0,1)$}
    \label{tab:tablerhypot}
    \begin{tabular}{l|l|l}
  & Naive & Compensated \\
      \hline
Zero ulp errors & 78.866  &  100\\
One ulp errors & 21.133 & 0
    \end{tabular}
  \end{center}
\end{table}

\section{Computing a Givens Rotation - {\tt DLARTG}}

In \cite{10.1145/567806.567809}, the authors describe a standard algorithm for computing real Givens rotations to zero out the second element in the vector $\left[ f , g \right]^T$. Note that we are using a slightly different element labeling than is normal because it will clarify our construction. The standard approach begins by dealing with the two exceptional cases where either $f$ or $g$ is zero and we will follow this convention. The next step involves corrections of scale to prevent avoidable floating point exceptions and there is some effort given to observe certain sign conventions, but in essence the approach they describe is to simply compute

\vspace{.2in}
\begin{algorithm} Naive DLARTG

\hrule
\begin{algorithmic}
   \State {$h = sqrt(f^2+g^2)$}
	\State {$c = f/h$}
	\State {$s = g/h$}
\end{algorithmic}
\hrule
\label{NaiveDLARTG}
\end{algorithm}
\vspace{.2in}

Note that the algorithm can be {\em mathematically} reduced to the problem of multiplying the elements by the reciprocal hypotenuse, and our compensated algorithm will work in precisely that manner. In particular, we will compute the reciprocal hypotenuse and its correction using the algorithm from the last section and then multiply these by the original $f$ and $g$ in a careful manner using the FMA. Note that the compensated version of the full DLARTG algorithm that appears below only differs from the rhypot code in the first two and last three lines.

\vspace{.2in}
\begin{algorithm} Compensated DLARTG

\hrule
\begin{algorithmic}
	\State {$x = abs(f)$}
	\State {$y = abs(g)$}
	\If{$x<y$} 
	\State {$swap(x,y)$}
   \EndIf
	\State {$x_{sq} = x^2$}
	\State {$y_{sq} = y^2$}
	\State {$\sigma = x_{sq}+y_{sq}$}
	\State {$\sigma_e = y_{sq}-(sigma-x_{sq})+fma(x,x,-x_{sq})+fma(y,y,-y_{sq})$}
   \State {$r = 1/\sigma$}
   \State {$\sigma = fma(-r,\sigma_e,fma(-r,\sigma,1))$}
	\State {$\rho = sqrt(\sigma)$}
	\State {$\tau = fma(-\rho,\rho,r)$}
	\State {$\bar{\nu} = \rho*fma(\sigma,\tau,\sigma)/2$}
   \State {$c = fma(f,\rho,f*\bar{\nu})$}
   \State {$s = fma(g,\rho,g*\bar{\nu})$}
\end{algorithmic}
\hrule
\label{CorrectedDLARTG}
\end{algorithm}
\vspace{.2in}

Once again, as a baseline for testing purposes we will first perform the naive calculation using the BigFloat format in Julia to compute $\bar{c}$ and $\bar{s}$. And once again we note that this computation cannot be guaranteed to represent a correctly rounded value of the true result (although it will do so in nearly every case). We will use $10^{9}$ normally distributed random inputs, that is both $x,y \sim \mathcal{N}(0,1)$. In the table below we summarize the error rates of each algorithm as before.

\begin{table}[h]
  \begin{center}
    \caption{Error rate for computing the Givens rotation with $x,y \sim \mathcal{N}(0,1)$}
    \label{tab:tableGivens}
    \begin{tabular}{l|l|l|l|l}
  &  \multicolumn{2}{c}{Naive} & \multicolumn{2}{c}{Compensated} \\
 & Cosine & Sine & Cosine & Sine \\
      \hline
Zero ulp errors & 66.563  & 66.567 & 100 & 100\\
One ulp errors & 33.207 &  33.204 & 0 & 0\\
Two ulp errors & 0.230 & 0.230 & 0 & 0 \\
    \end{tabular}
  \end{center}
\end{table}

\section{Conclusions}

We have shown that it is possible to compute the exact value of the standard Newton step for reciprocal square root for estimates $y$ that are sufficiently close to the true value. This Newton step can be used to get an estimate that is weakly rounded which is a strong guarantee of high accuracy. We then showed how the exact Newton step can be modified to get a better estimate of the actual error of $y$. This leads to an algorithm that we conjecture is correctly rounded although we cannot prove it at this time. Moreover, this modified compensation can be shown to significantly improve the accuracy of certain square-root free methods for estimating the reciprocal square root. Finally, we saw that these methods can be applied to the calculation of the reciprocal hypotenuse and Givens rotations to generate values that are highly accurate experimentally. 

\section{Acknowledgements}
The author wishes to sincerely thank the reviewers who made many meaningful and insightful suggestions that greatly improved this work.

\bibliography{hyp}{}
\bibliographystyle{acm}

\end{document}